\def\latexversion#1#2{#2} 
\latexversion{\def\smathsf{\mathsf}} 
{\font\sf=cmss10\def\mathsf#1{{\mbox{\sf#1}}}
\font\smsf=cmss8\def\smathsf#1{{\mbox{\smsf#1}}}} 

\documentstyle[twoside]{article}

\catcode`\@=11
\long\def\@makefntext#1{
\protect\noindent \hbox to 3.2pt {\hskip-.9pt
$^{{\eightrm\@thefnmark}}$\hfil}#1\hfill}		

\def\thefootnote{\fnsymbol{footnote}}
\def\@makefnmark{\hbox to 0pt{$^{\@thefnmark}$\hss}}	

\def\ps@myheadings{\let\@mkboth\@gobbletwo
\def\@oddhead{\hbox{}
\rightmark\hfil\eightrm\thepage}
\def\@oddfoot{}\def\@evenhead{\eightrm\thepage\hfil
\leftmark\hbox{}}\def\@evenfoot{}
\def\sectionmark##1{}\def\subsectionmark##1{}}

\oddsidemargin=\evensidemargin
\renewcommand{\thefootnote}{\fnsymbol{footnote}}

\newcounter{sectionc}\newcounter{subsectionc}\newcounter{subsubsectionc}
\renewcommand{\section}[1] {\vspace{12pt}\addtocounter{sectionc}{1}
\setcounter{subsectionc}{0}\setcounter{subsubsectionc}{0}\noindent
	{\tenbf\thesectionc. #1}\par\vspace{5pt}}
\renewcommand{\subsection}[1] {\vspace{12pt}\addtocounter{subsectionc}{1}
	\setcounter{subsubsectionc}{0}\noindent
	{\bf\thesectionc.\thesubsectionc. {\kern1pt \bfit #1}}\par\vspace{5pt}}
\renewcommand{\subsubsection}[1] {\vspace{12pt}\addtocounter{subsubsectionc}{1}
	\noindent{\tenrm\thesectionc.\thesubsectionc.\thesubsubsectionc.
	{\kern1pt \tenit #1}}\par\vspace{5pt}}
\newcommand{\nonumsection}[1] {\vspace{12pt}\noindent{\tenbf #1}
	\par\vspace{5pt}}

\newcounter{appendixc}
\newcounter{subappendixc}[appendixc]
\newcounter{subsubappendixc}[subappendixc]
\renewcommand{\thesubappendixc}{\Alph{appendixc}.\arabic{subappendixc}}
\renewcommand{\thesubsubappendixc}
	{\Alph{appendixc}.\arabic{subappendixc}.\arabic{subsubappendixc}}

\renewcommand{\appendix}[1] {\vspace{12pt}
        \refstepcounter{appendixc}
        \setcounter{figure}{0}
        \setcounter{table}{0}
        \setcounter{lemma}{0}
        \setcounter{theorem}{0}
        \setcounter{corollary}{0}
        \setcounter{definition}{0}
        \setcounter{equation}{0}
        \renewcommand{\thefigure}{\Alph{appendixc}.\arabic{figure}}
        \renewcommand{\thetable}{\Alph{appendixc}.\arabic{table}}
        \renewcommand{\theappendixc}{\Alph{appendixc}}
        \renewcommand{\thelemma}{\Alph{appendixc}.\arabic{lemma}}
        \renewcommand{\thetheorem}{\Alph{appendixc}.\arabic{theorem}}
        \renewcommand{\thedefinition}{\Alph{appendixc}.\arabic{definition}}
        \renewcommand{\thecorollary}{\Alph{appendixc}.\arabic{corollary}}
        \renewcommand{\theequation}{\Alph{appendixc}.\arabic{equation}}
        \noindent{\tenbf Appendix \theappendixc #1}\par\vspace{5pt}}
\newcommand{\subappendix}[1] {\vspace{12pt}
        \refstepcounter{subappendixc}
        \noindent{\bf Appendix \thesubappendixc. {\kern1pt \bfit #1}}
	\par\vspace{5pt}}
\newcommand{\subsubappendix}[1] {\vspace{12pt}
        \refstepcounter{subsubappendixc}
        \noindent{\rm Appendix \thesubsubappendixc. {\kern1pt \tenit #1}}
	\par\vspace{5pt}}

\newcommand{\textlineskip}{\baselineskip=13pt}
\newcommand{\smalllineskip}{\baselineskip=10pt}

\def\eightcirc{
\begin{picture}(0,0)
\put(4.4,1.8){\circle{6.5}}
\end{picture}}
\def\eightcopyright{\eightcirc\kern2.7pt\hbox{\eightrm c}}

\newcommand{\copyrightheading}[1]
	{\vspace*{-2.5cm}\smalllineskip{\flushleft
	{\footnotesize International Journal of Modern Physics B, #1}\\
	{\footnotesize $\eightcopyright$\, World Scientific Publishing
	 Company}\\
	 }}

\newcommand{\pub}[1]{{\begin{center}\footnotesize\smalllineskip
	Received #1\\
	\end{center}
	}}

\def\abstracts#1#2#3{{
	\centering{\begin{minipage}{4.5in}\baselineskip=10pt\footnotesize
	\parindent=0pt #1\par
	\parindent=15pt #2\par
	\parindent=15pt #3
	\end{minipage}}\par}}

\def\keywords#1{{
	\centering{\begin{minipage}{4.5in}\baselineskip=10pt\footnotesize
	{\footnotesize\it Keywords}\/: #1
	\end{minipage}}\par}}

\renewenvironment{thebibliography}[1]			
	{\frenchspacing
	 \ninerm\baselineskip=11pt
	 \begin{list}{\arabic{enumi}.}
	{\usecounter{enumi}\setlength{\parsep}{0pt}
	 \setlength{\leftmargin 12.7pt}{\rightmargin 0pt} 
	 \setlength{\itemsep}{0pt} \settowidth
	{\labelwidth}{#1.}\sloppy}}{\end{list}}

\newcounter{itemlistc}
\newcounter{romanlistc}
\newcounter{alphlistc}
\newcounter{arabiclistc}

\newcommand{\fcaption}[1]{
        \refstepcounter{figure}
        \setbox\@tempboxa = \hbox{\footnotesize Fig.~\thefigure. #1}
        \ifdim \wd\@tempboxa > 5in
           {\begin{center}
        \parbox{5in}{\footnotesize\smalllineskip Fig.~\thefigure. #1}
            \end{center}}
        \else
             {\begin{center}
             {\footnotesize Fig.~\thefigure. #1}
              \end{center}}
        \fi}

\newcommand{\tcaption}[1]{
        \refstepcounter{table}
        \setbox\@tempboxa = \hbox{\footnotesize Table~\thetable. #1}
        \ifdim \wd\@tempboxa > 5in
           {\begin{center}
        \parbox{5in}{\footnotesize\smalllineskip Table~\thetable. #1}
            \end{center}}
        \else
             {\begin{center}
             {\footnotesize Table~\thetable. #1}
              \end{center}}
        \fi}

\def\@citex[#1]#2{\if@filesw\immediate\write\@auxout
	{\string\citation{#2}}\fi
\def\@citea{}\@cite{\@for\@citeb:=#2\do
	{\@citea\def\@citea{,}\@ifundefined
	{b@\@citeb}{{\bf ?}\@warning
	{Citation `\@citeb' on page \thepage \space undefined}}
	{\csname b@\@citeb\endcsname}}}{#1}}

\newif\if@cghi
\def\cite{\@cghitrue\@ifnextchar [{\@tempswatrue
	\@citex}{\@tempswafalse\@citex[]}}
\def\citelow{\@cghifalse\@ifnextchar [{\@tempswatrue
	\@citex}{\@tempswafalse\@citex[]}}
\def\@cite#1#2{{$\null^{#1}$\if@tempswa\typeout
	{IJCGA warning: optional citation argument
	ignored: `#2'} \fi}}

\def\pmb#1{\setbox0=\hbox{#1}
	\kern-.025em\copy0\kern-\wd0
	\kern.05em\copy0\kern-\wd0
	\kern-.025em\raise.0433em\box0}

\def\fnt#1#2{\footnotetext{\kern-.3em
	{$^{\mbox{\scriptsize #1}}$}{#2}}}

\def\fpage#1{\begingroup
\voffset=.3in
\thispagestyle{empty}\begin{table}[b]\centerline{\footnotesize #1}
	\end{table}\endgroup}

\def\runninghead#1#2{\pagestyle{myheadings}
\markboth{{\protect\footnotesize\it{\quad #1}}\hfill}
{\hfill{\protect\footnotesize\it{#2\quad}}}}
\font\tenrm=cmr10
\font\tenit=cmti10
\font\tenbf=cmbx10
\font\bfit=cmbxti10 at 10pt
\font\ninerm=cmr9
\font\nineit=cmti9
\font\ninebf=cmbx9
\font\eightrm=cmr8

\def\qed{\hbox{${\vcenter{\vbox{			
   \hrule height 0.4pt\hbox{\vrule width 0.4pt height 6pt
   \kern5pt\vrule width 0.4pt}\hrule height 0.4pt}}}$}}

\renewcommand{\thefootnote}{\fnsymbol{footnote}}	

\def\bsc{{\sc a\kern-6.4pt\sc a\kern-6.4pt\sc a}}	
\def\bflatex{\bf L\kern-.30em\raise.3ex\hbox{\bsc}\kern-.14em
T\kern-.1667em\lower.7ex\hbox{E}\kern-.125em X}

\def\Journal#1#2#3#4{{#1}\ {\ninebf #2}, #3 (#4)}
\def\myem{\nineit}
\newcount\mya\newcount\myb\font\my=cmr10 at 15pt
\def\myphi{\hbox{\my\char'010}}
\def\hypp#1#2#3#4{\mya=#1\myb=#1\advance\myb by1{}_{\the\myb}
\myphi_\the\mya\left[{\omega,#2\atop\phantom{\omega,}#3};#4\right]}
\def\hypg#1#2#3{{}_{p+1}\myphi_{p}
\left[{#1\atop{#2\phantom{w}}};#3\right]}
\def\hypq#1#2#3#4{\mya=#1\myb=#1\advance\myb by1
{}_{\the\myb}\myphi_\the\mya\left[{#2\atop#3};#4\right]}
\def\lbar#1{\overline{#1}}
\def\wb{\lbar{W}}
\def\betaE#1{\displaystyle{#1\over k_{\rm B}T}}
\def\half{{\textstyle{1\over2}}}
\def\halfs{{\scriptstyle\frac12}}
\newcommand{\mysubsection}{\subsection}

\begin{document}

\input{psfig}

\runninghead{H.\ Au-Yang \&\ J.\ H.\ H. Perk}
{The Many Faces of the Chiral Potts Model}

\normalsize\textlineskip
\thispagestyle{empty}
\setcounter{page}{1}

\copyrightheading{}			

\vspace*{0.88truein}

\fpage{1}
\centerline{\bf THE MANY FACES OF THE CHIRAL POTTS MODEL}
\vspace*{0.37truein}
\centerline{\footnotesize HELEN AU-YANG and
JACQUES H.H.\ PERK\footnote{E-mail: perk@okstate.edu}}
\vspace*{0.015truein}
\centerline{\footnotesize\it Department of Physics,
Oklahoma State University}
\baselineskip=10pt
\centerline{\footnotesize\it 145 Physical Sciences, Stillwater,
OK 74078-3072, USA}
\vspace*{0.225truein}
\pub{26 August 1996}

\vspace*{0.21truein}
\abstracts{In this talk, we give a brief overview of
several aspects of the theory of the chiral Potts model, including
higher-genus solutions of the star--triangle and tetrahedron
equations, cyclic representations of affine quantum groups, basic
hypergeometric functions at root of unity, and possible
applications.}{}{}

\vspace*{10pt}


\vspace*{1pt}\textlineskip	
\section{Introduction}\label{sec:intro}
\vspace*{-0.5pt}
\noindent
The chiral Potts model is a generalization of the $N$-state Potts model
allowing for handedness of the pair interactions. While known under
several other names, it has received much interest in the past two
decades. It may have appeared first in a short paper by Wu and
Wang\cite{WuWang} in the context of duality transformations.

However, active studies of chiral Potts models did not start until a
few years after that, when \"Ostlund\cite{Os} and Huse\cite{Hu}
introduced the more special chiral clock model as a model for
commensurate--incommensurate phase transitions. Much has been published
since and in this talk we can only highlight some of the developments
and discuss a few of those related to the integrable manifold in more
detail.

\mysubsection{Domain wall theory of incommensurate states in adsorbed
layers}
\noindent
Following the original papers of \"Ostlund\cite{Os} and Huse\cite{Hu}
there was immediately much interest in their model, because it can be
used to describe wetting phenomena in commensurate phases, transitions
to incommensurate states and it provides through the domain wall theory
a model for adsorbed monolayers.\cite{HSF}$^-$\cite{AP-O}

\mysubsection{Chiral NN interactions to describe further neighbor
interactions}
\noindent
One may object that next-nearest and further neighbor interactions are
to be seen as the physical cause of incommensurate states rather than
chiral interactions.\break However, one can show that by a block-spin
transformation such longer-range models can be mapped to a
nearest-neighbor interaction model but with chiral
interactions\cite{Barber}$^-$\cite{AP-F} in general. From the viewpoint
of integrability, for example, the nearest-neighbor picture is
preferred.

\mysubsection{New solutions of quantum Lax pairs, or star--triangle
equations}
\noindent
The integrable chiral Potts model\cite{AMPTY}$^-$\cite{AP-tani}
provides new solutions of the Yang--Baxter equation, which could help
understanding chiral field directions and correlations in lattice and
conformal field theories. In fact, the belief that to each conformal
field theory corresponds a lattice model was the hidden motivation
behind the original discovery.\cite{AMPTY}

\mysubsection{Higher-genus spectral parameters (Fermat curves)}
\noindent
The integrable chiral Potts models are different from all other
solvable models based on the Yang--Baxter equations. The spectral
parameters (or rapidities) lie on higher-genus
curves.\cite{AMPTY}$^-$\cite{AP-tani}

\mysubsection{Level crossings in quantum chains (1D vs 2D)}
\noindent
The physics of incommensurate states in one-dimensional quantum chiral
Potts chains is driven by level crossings,\cite{AMP} which are
forbidden in the classical case by\break the Perron--Frobenius theorem.
The two-dimensional chiral Potts model has its\break integrable
submanifold within the commensurate phase, which ends at the
Fateev--Zamolodchikov multicritical point.\cite{FZ}

\mysubsection{Exact solutions for several physical quantities}
\noindent
All other Yang--Baxter solvable models discovered so far have a
uniformization that is  based on elementary or elliptic functions,
with meromorphic dependences on differences (and sums) of spectral
parameters. This is instrumental in the\break evaluation of their
physical quantities. For the integrable chiral Potts model\break
several quantities have been obtained using new approaches without
using an\break explicit
uniformization.\cite{AMP,AMPT}$^-$\cite{oRB,AP-O,AP-tani}

\mysubsection{Multi-component versions}
\noindent
Multicomponent versions of the chiral Potts model\cite{AP-multi} may
be of interest in many fields of study, such as the structure of lipid
bilayers for which the Pearce--Scott model has been
introduced.\cite{MPS}

\mysubsection{Large N-limit}
\noindent
The integrable chiral Potts model allows three large-$N$ limits that
may be useful in connection with $W_{\infty}$ algebras.\cite{AP-inf}

\mysubsection{Cyclic representations of quantum groups at roots of 1}
\noindent
The chiral Potts model can be viewed as the first application of the
theory of cyclic representations\cite{dCK} of affine quantum groups.


\textheight=7.8truein
\setcounter{footnote}{0}
\renewcommand{\thefootnote}{\alph{footnote}}
\mysubsection{Generalizing free fermions (Onsager algebra)}
\noindent
The operators in the superintegrable subcase of the chiral Potts model
obey\break Onsager's loop group algebra, making the model integrable
for the same two reasons as the $N=2$ Ising or free-fermion model.
Howes, Kadanoff, and den Nijs\cite{HKdN} first noted special features of
series expansions for a special case of the $N=3$ quantum chain. This
was generalized by von Gehlen and Rittenberg\cite{vGR} to arbitrary $N$
using the Dolan--Grady criterium, which was later explained to be
equivalent to Onsager's loop algebra relations.\cite{PD,Dav}

\mysubsection{Solutions of tetrahedron equations}
\noindent
Bazhanov and Baxter\cite{BB}$^-$\cite{SMS} have shown that the
sl($n$) generalization\cite{DJMM}$^-$\cite{KM} of the integrable
chiral Potts model can be viewed as an $n$-layered $N$-state
generalization of the three-dimensional Zamolodchikov model.

\mysubsection{Cyclic hypergeometric functions}
\noindent
Related is a new theory of basic hypergeometric functions at root
of unity discussed in some detail at the end of this talk.

\mysubsection{New models with few parameters}
\noindent
In the integrable submanifold with positive Boltzmann weights several
ratios of the parameters are nearly constant, suggesting the study of
new two-parameter $N$-state models.\cite{AP-O,AP-C}

\section{The Chiral Potts Model}
\begin{figure}[htbp]
\hbox{\hspace*{.17in}\vbox{\vspace*{.15in}%
\psfig{figure=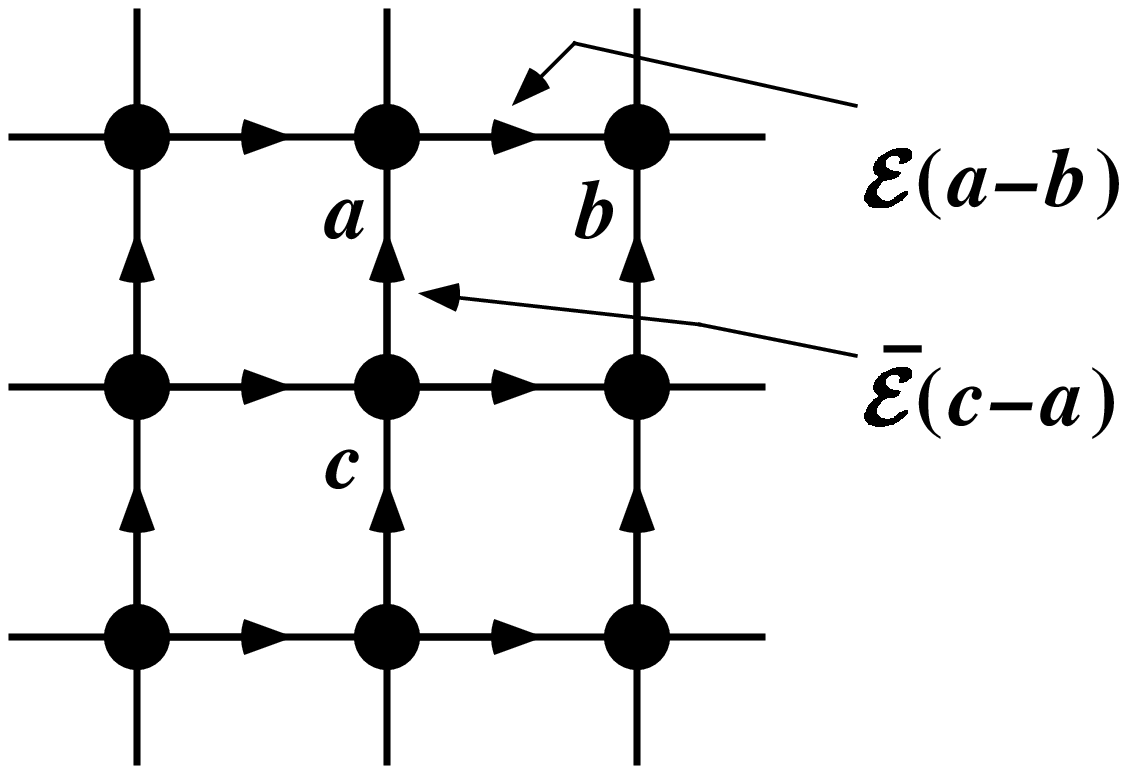,height=1.1in}%
\vspace*{.35in}}\hspace*{.35in}%
\psfig{figure=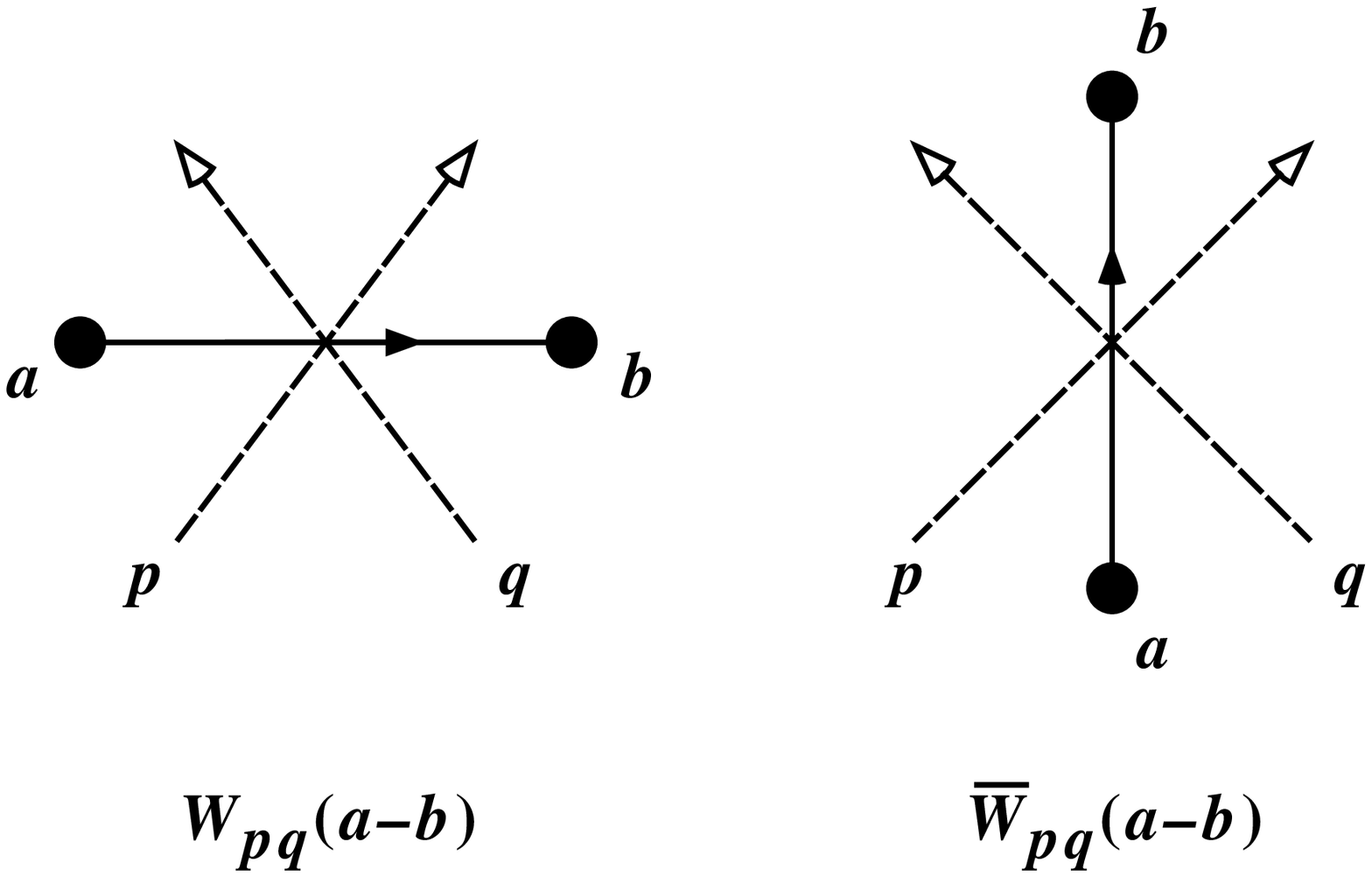,height=1.7in}}
\vspace*{13pt}
\fcaption{Interaction energies and Boltzmann weights of horizontal and
vertical couplings in the chiral Potts model. The subtraction of the
state variables $a,b,c=1,2,\ldots,N$ at the lattice sites is done
modulo $N$. The dashed lines with open arrows denote oriented lines on
the medial lattice where the rapidity variables $p,q,\ldots,$ of the
integrable submanifold live.\hfill
\label{fig2}}
\end{figure}

\noindent
The most general chiral Potts model is defined on a graph or lattice,
see Fig.~\ref{fig2}, where the interaction energies
\begin{equation}
{\cal E}(n)={\cal E}(n+N)=\sum_{j=1}^{N-1}E_j\,\omega^{jn},
\qquad \omega\equiv e^{2\pi i/N},
\label{en-cp}
\end{equation}
depend on the differences $n=a-b$ mod$\,N$ of the spin variables
$a$ and $b$ on two neighboring sites. We can write
\begin{equation}
\betaE{E^{\vphantom{\displaystyle\ast}}_j}=
\betaE{E^{\displaystyle\ast}_{N-j}}=-K_j\,\omega^{\Delta_j},
\quad j=1,\ldots,\lfloor\half N\rfloor,
\label{en-cp2}
\end{equation}
where $K_j$ and $\Delta_j$ constitute $N-1$ independent variables.
Then, for $N$ odd we have a sum of ``clock model" terms
\begin{equation}
-\,\betaE{{\cal E}(n)}=
\sum_{j=1}^{\halfs(N-1)}2K_j
\cos\Big[{2\pi\over N}(jn+\Delta_j)\Big],
\label{en-odd}
\end{equation}
whereas, for $N$ even we have an additional Ising term, {\it i.e.}
\begin{equation}
-\,\betaE{{\cal E}(n)}=\sum_{j=1}^{\halfs N-1}2K_j
\cos\Big[{2\pi\over N}(jn+\Delta_j)\Big]+K_{\halfs N}(-1)^n.
\label{en-even}
\end{equation}
The Boltzmann weight corresponding to the edge is given by
\begin{equation}
W(n)=e^{-{\cal E}(n)/k_{\rm B}T}_{\vphantom{X}}.
\label{en-bw}
\end{equation}

\section{The Integrable Chiral Potts Model}
\begin{figure}[htbp]
\hbox{\hspace*{.15in}\psfig{figure=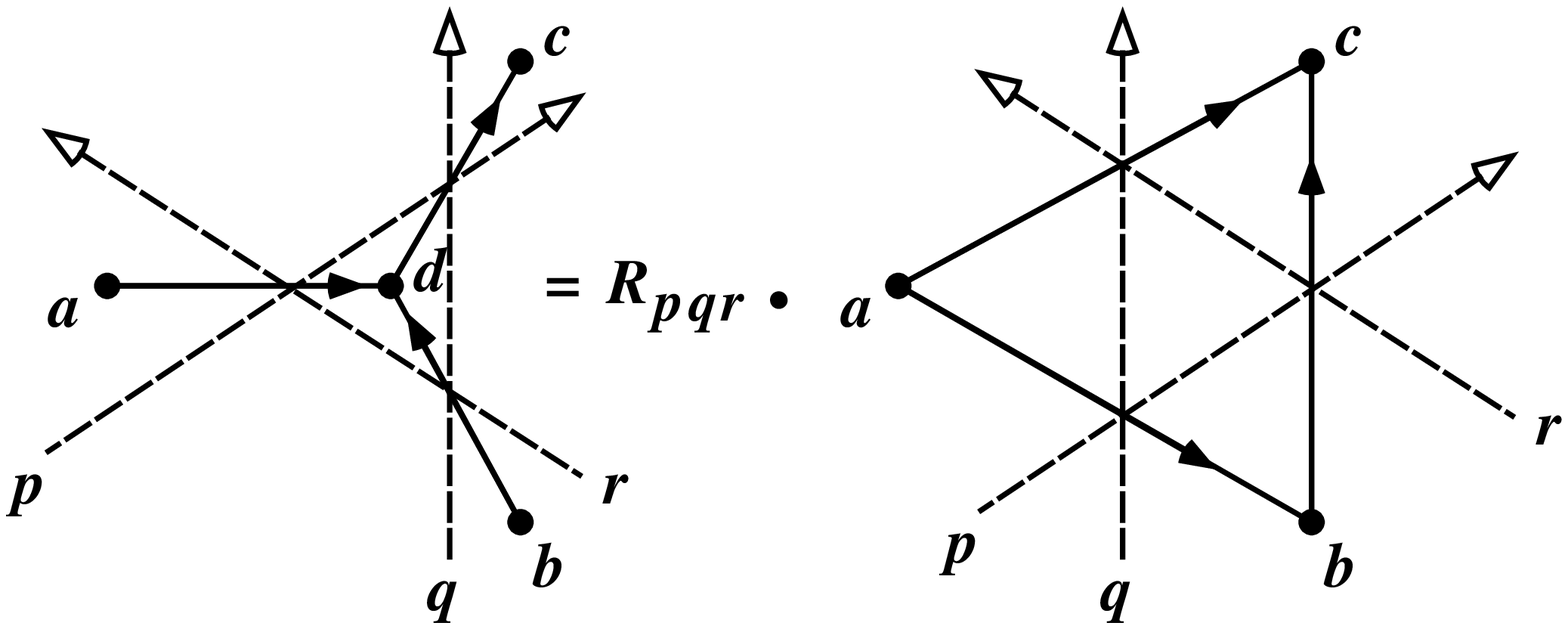,height=1.9in}}
\vspace*{13pt}
\fcaption{Star--Triangle Equation.
\label{fig3}}
\end{figure}

\noindent
In the integrable chiral Potts model, we have besides ``spins"
$a,b,\ldots$ defined mod$\,N$, ``rapidity lines" $p,q,\ldots$
all pointing in one halfplane. The weights satisfy the star--triangle
equation, see also Fig.~\ref{fig3},
\begin{equation}
\sum^{N}_{d=1}\wb_{qr}(b-d)W_{pr}(a-d)\wb_{pq}(d-c)
=R_{pqr}W_{pq}(a-b)\wb_{pr}(b-c)W_{qr}(a-c).
\label{STE}
\end{equation}
In full generality there are six sets of weights to be found and a
constant $R$. But the solution is in terms of two functions depending
on spin differences and pairs of rapidity variables:
\begin{eqnarray}
W_{pq}(n)&=&W_{pq}(0)\prod^{n}_{j=1}\biggl({\mu_p\over\mu_q}\cdot
{y_q-x_p\omega^j\over y_p-x_q\omega^j}\biggr),\nonumber\\
\wb_{pq}(n)&=&\wb_{pq}(0)\prod^{n}_{j=1}\biggl(\mu_p\mu_q\cdot
{\omega x_p-x_q\omega^j\over y_q-y_p\omega^j}\biggr),
\label{weights}
\end{eqnarray}
with $R$ depending on three rapidity variables.
Periodicity modulo $N$ gives for all rapidity pairs $p$ and $q$
\begin{equation}
\biggl({\mu_p\over\mu_q}\biggr)^N={y_p^N-x_q^N\over y_q^N-x_p^N},
\qquad(\mu_p\mu_q)^N={y_q^N-y_p^N\over x_p^N-x_q^N}.
\label{periodic}
\end{equation}
Hence, we can define (rescale) $k,k'$ such that
\begin{equation}
\mu_p^N={k'\over1-k\,x_p^N}={1-k\,y_p^N\over k'},\quad
x_p^N+y_p^N=k(1+x_p^Ny_p^N),\quad k^2+k'^2=1,
\label{curve}
\end{equation}
and similarly with $p$ replaced by $q$. The rapidities live on
a curve of genus $g>1$ that is of Fermat type.

\section{Physical Cases}
\noindent
There are two conditions for physical cases:
\renewcommand{\labelenumi}{\Roman{enumi}.}
\begin{enumerate}
\item{Planar Model with Real Positive Boltzmann Weights,}
\item{Hermitian Quantum Spin Chain.}
\end{enumerate}
Usually, as in the six-vertex model in electric fields where the
quantum chain would have either imaginary or real Dzyaloshinsky--Moriya
interactions, one cannot require both simultaneously. Only for the
nonchiral (reflection-positive) subcase of the Fateev--Zamolodchikov
model\cite{FZ} are both physical conditions simultaneously fulfilled
for $N>2$.

We note that the Hermitian quantum spin chain submanifold contains
the\break superintegrable case, where both the star--triangle (or
Yang--Baxter) equation and the Onsager algebra are satisfied.

\section{Generalization of Free Fermion Model}
\begin{figure}[htbp]
\hbox{\hspace*{0.92in}\psfig{figure=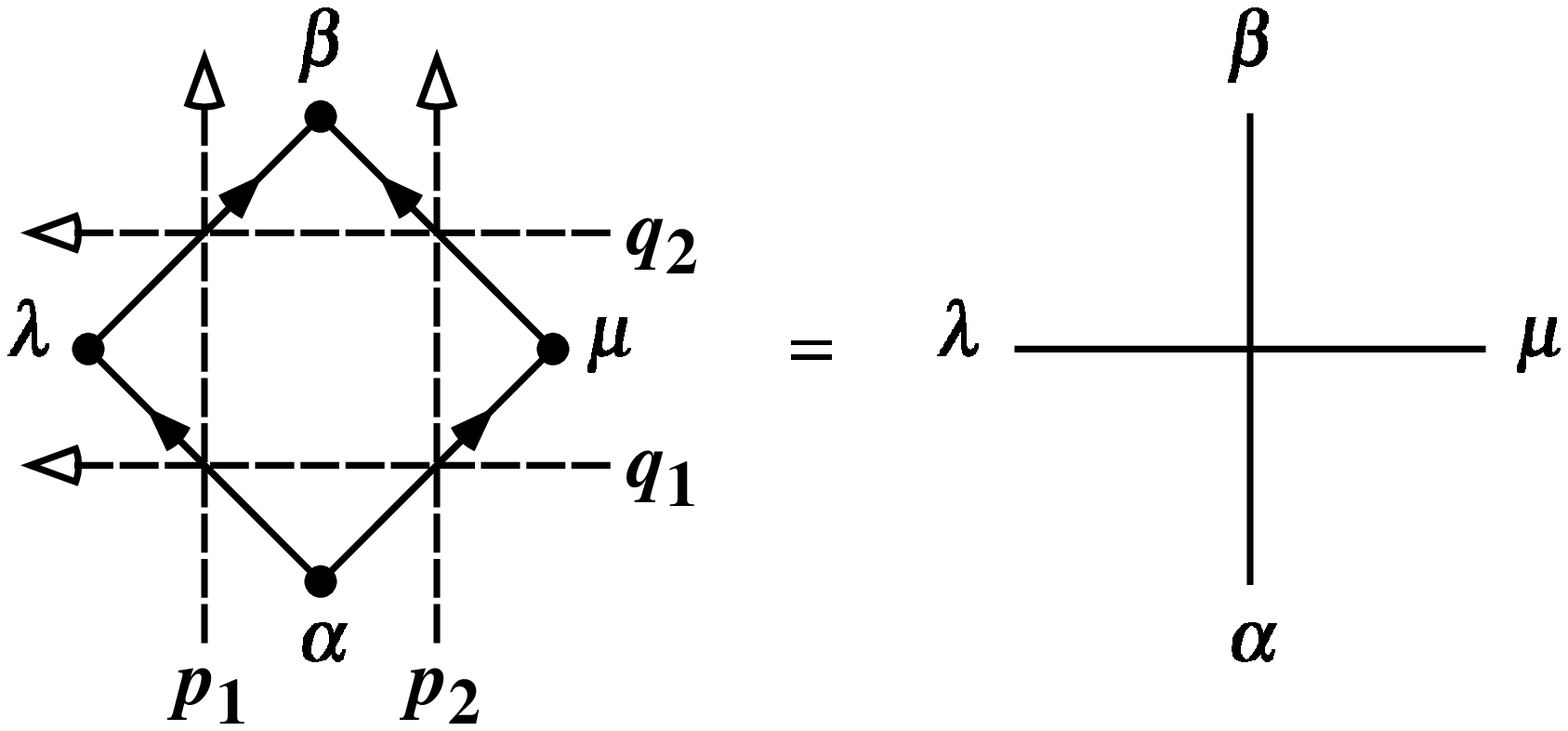,height=1.5in}}
\vspace*{13pt}
\fcaption{Vertex model satisfying the Yang--Baxter Equation.
\label{fig4}}
\end{figure}

\noindent
Combining four chiral Potts model weights in a square, as in
Fig.~\ref{fig4},
we obtain the $\mathsf{R}$-matrix of the $N$-state
generalization\cite{BPA} of the checkerboard Ising model
\begin{equation}
R_{\alpha\beta|\lambda\mu}=
\wb_{p_1q_1}(\alpha-\lambda)\wb_{p_2q_2}(\mu-\beta)
W_{p_2q_1}(\alpha-\mu)W_{p_1q_2}(\lambda-\beta).
\label{ff-a}
\end{equation}
Applying a Fourier transform gauge transformation we obtain\cite{BPA}
from this
\begin{equation}
{\hat R}_{\alpha\beta|\lambda\mu}=
{s_\beta\,t_\lambda\over s_\alpha\,t_\mu}\,{1\over N^2}\,
\sum^{N}_{\alpha'=1}\sum^{N}_{\beta'=1}
\sum^{N}_{\lambda'=1}\sum^{N}_{\mu'=1}
\omega^{-\alpha\alpha'+\beta\beta'+\lambda\lambda'-\mu\mu'}
R_{\alpha'\beta'|\lambda'\mu'},
\label{ff-b}
\end{equation}
which is the $\mathsf{R}$-matrix of an $N$-state generalization of the
free-fermion eight-vertex model for $N=2$. Indeed,
${\hat R}_{\alpha\beta|\lambda\mu}$ is nonzero only if
$\lambda+\beta=\alpha+\mu$ mod$\,N$ which generalizes the eight-vertex
condition.

In eq.~(\ref{ff-b}) $s_{\alpha}$ and $t_{\lambda}$ are free parameters,
corresponding to gauge freedom, which may be edge dependent.

\section{Transfer Matrices}
\noindent
\begin{figure}[htbp]
\hbox{\hspace*{0.39in}\psfig{figure=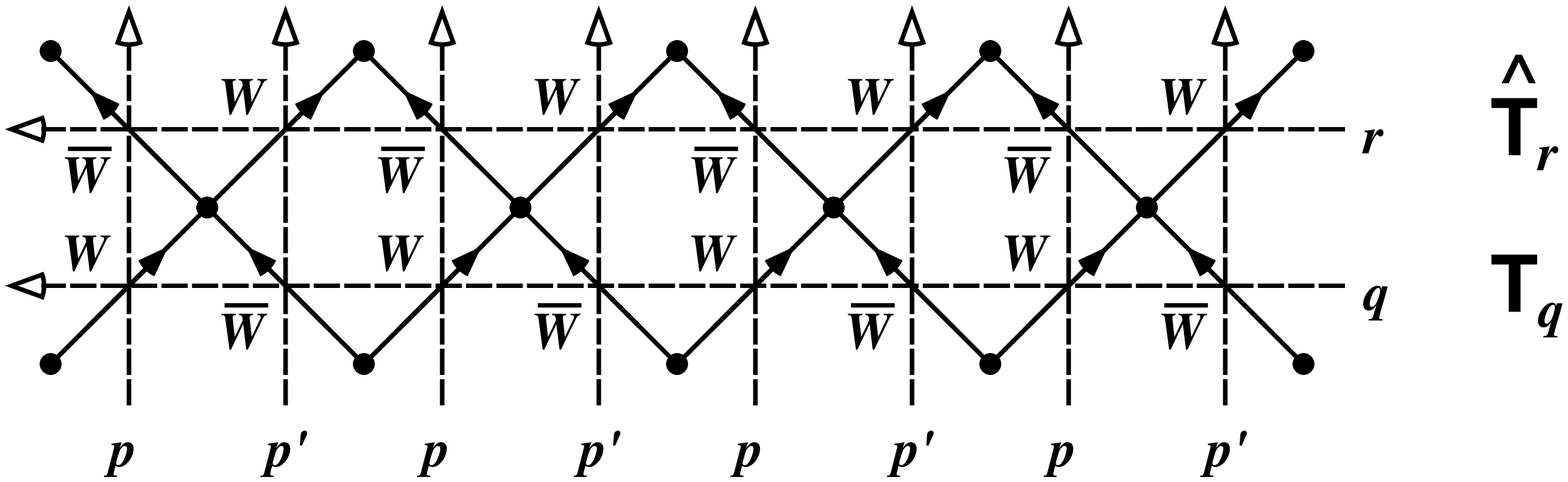,height=1.3in}}
\vspace*{13pt}
\fcaption{Diagonal-to-diagonal transfer matrices
$\smathsf{T}_q$ and $\hat{\smathsf{T}}_r$.
\label{fig5}}
\end{figure}
If the weight functions $W$ and $\wb$ are solutions of the star--triangle
equation (\ref{STE}) and parametrized as in (\ref{weights}) with
variables on rapidity lines lying on the algebraic curve (\ref{curve}), all
diagonal transfer matrices commute. This is depicted in Fig.~\ref{fig5},
where the vertical rapidity variables have been chosen alternatingly;
more generally they could all be chosen independently.

\begin{figure}[htbp]
\hbox{\hspace*{0.17in}\psfig{figure=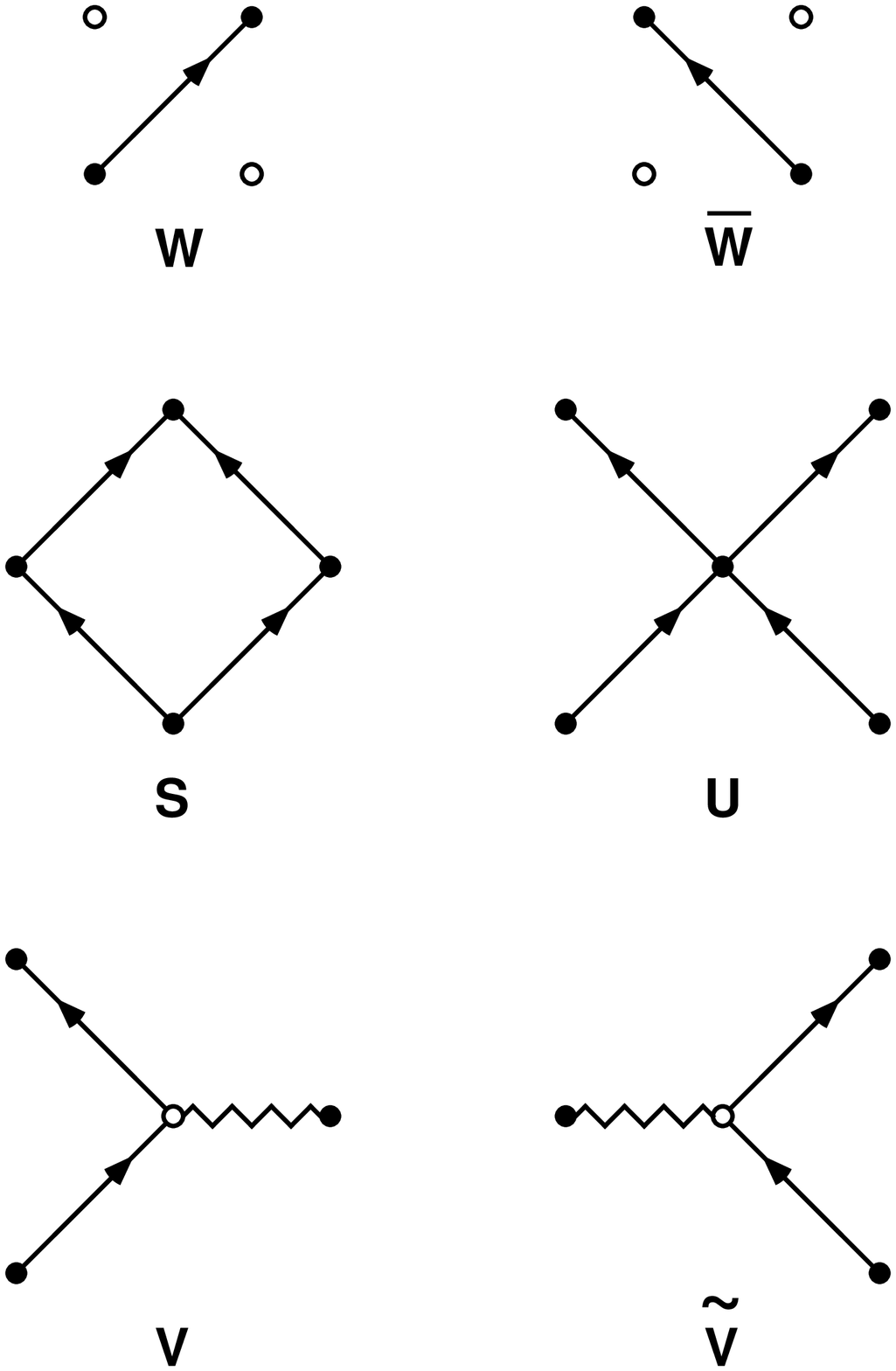,height=2.2in}%
\hspace*{.6in}%
\psfig{figure=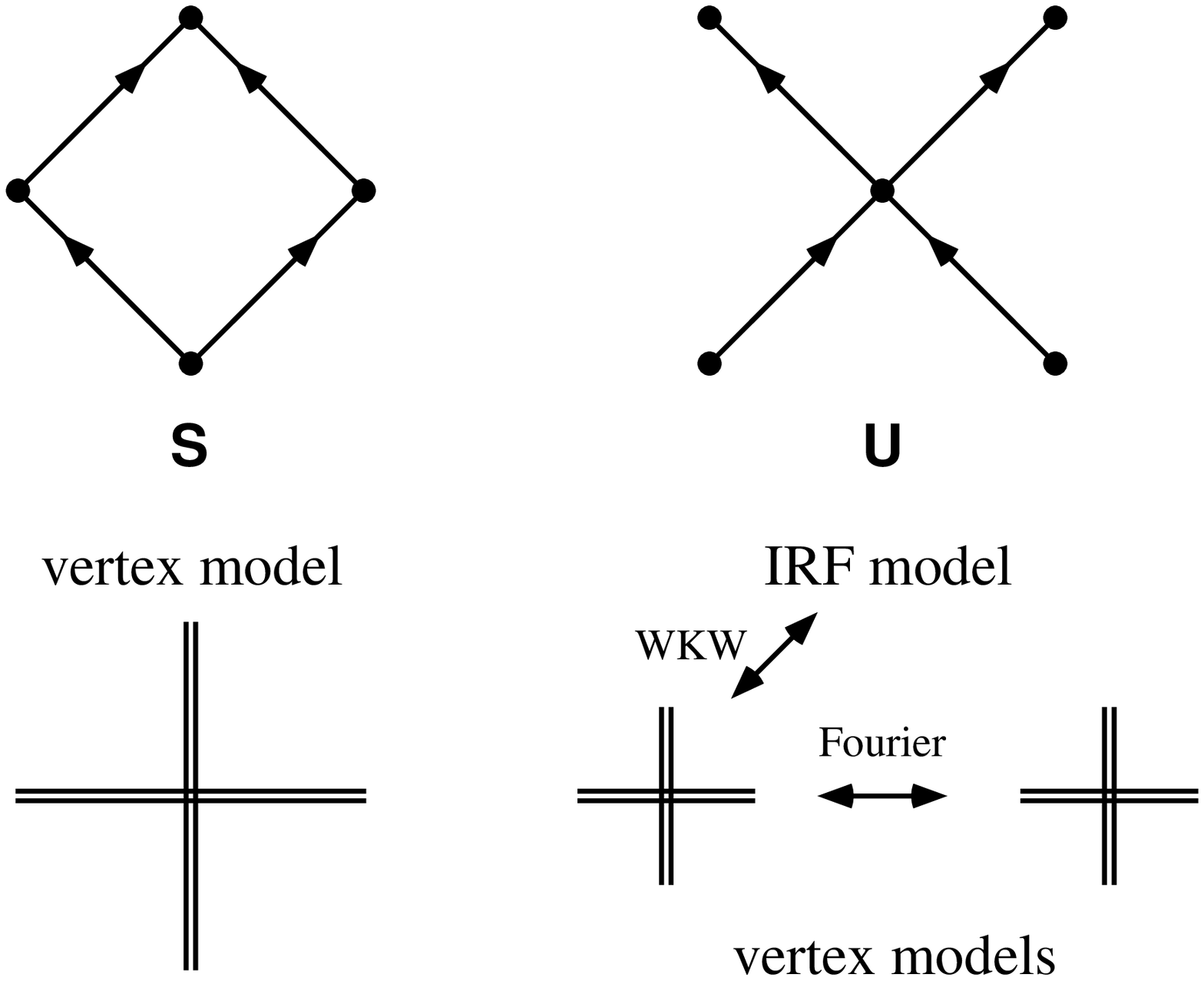,height=2.2in}}
\vspace*{13pt}
\fcaption{Various $\smathsf{R}$-matrices related to the chiral Potts
model. Vertex weight $\smathsf{S}$ and IRF weight $\smathsf{U}$ are
built from two $W$'s and two $\wb$'s and have two horizontal and two
vertical rapidity lines. By the Wu--Kadanoff-Wegner map, IRF weight
$\smathsf{U}$ is also a vertex weight. The two zigzag lines in
$\smathsf{V}$ and $\tilde{\smathsf{V}}$ represent Fourier and inverse
Fourier transform, respectively.
\label{fig7}}
\end{figure}
We note that there are several ways to introduce
$\mathsf{R}$-matrices, as is indicated in Fig.~\ref{fig7}. First, we
have the original weights $\mathsf{W}$ and $\lbar{\mathsf{W}}$. Next,
we have the vertex weight (square) $\mathsf{S}$ and the IRF weight
(star) $\mathsf{U}$, both consisting of two $W$ and two $\wb$ weights.
Finally, we have the three-spin interactions $\mathsf{V}$ and
$\tilde{\mathsf{V}}$, its transposed or rotated version, which play a
special role in the theory.\cite{AMPTY,AP-tani} Note that $\mathsf{S}$
and $\mathsf{U}$ can be constructed from a $\mathsf{V}$ and a
$\tilde{\mathsf{V}}$, for example as
$\mathsf{U}=\mathsf{V}\cdot\tilde{\mathsf{V}}$.

\section{The Construction of Bazhanov and Stroganov}
\begin{figure}[htbp]
\hbox{\hspace*{.3in}\psfig{figure=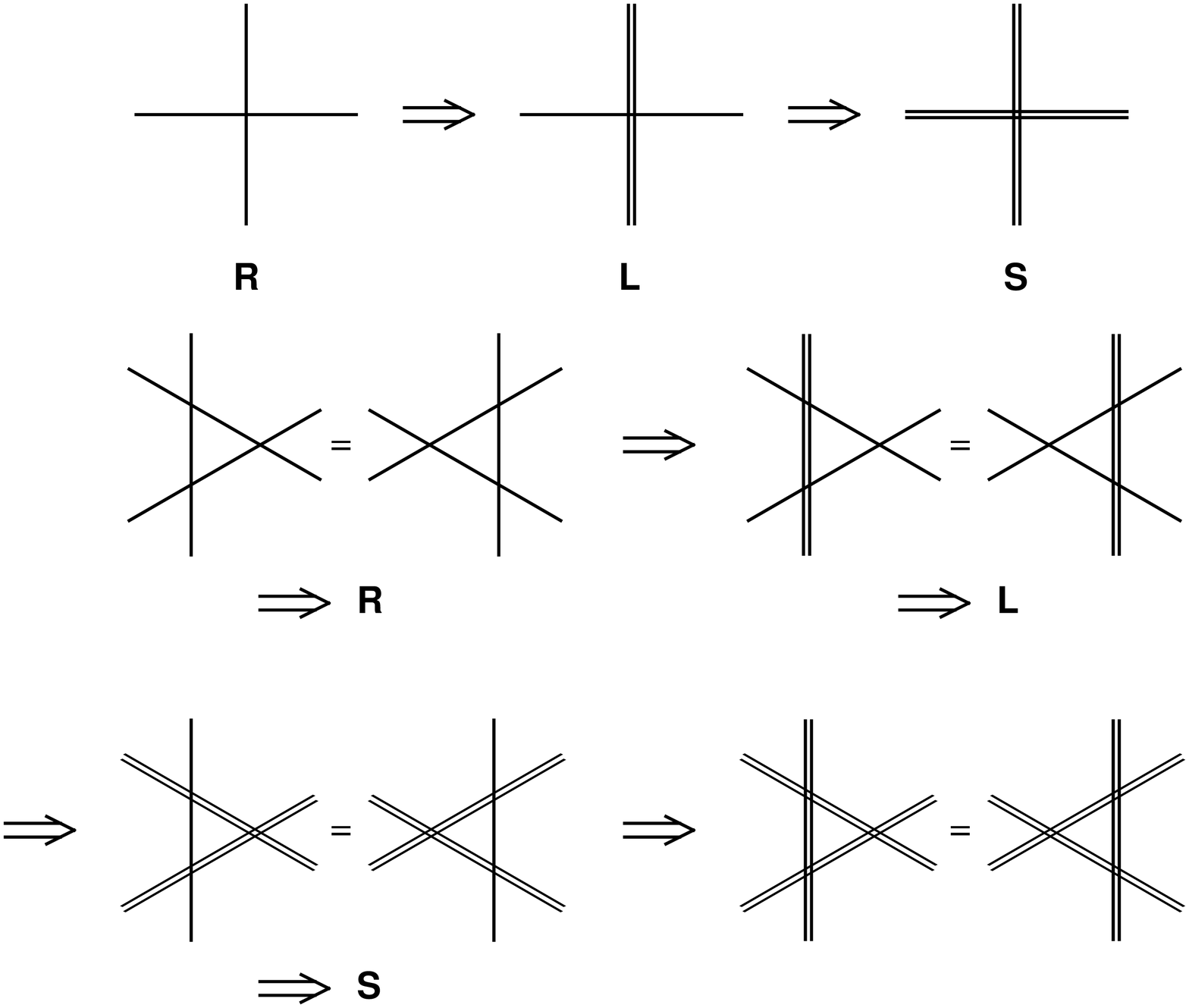,height=3.5in}}
\vspace*{13pt}
\fcaption{The construction of Bazhanov and Stroganov. Single lines
correspond to spin $\halfs$ highest-weight representations; double
lines correspond to minimal cyclic representations.
\label{fig8}}
\end{figure}

\noindent
Bazhanov and Stroganov\cite{BS} have shown that the
$\mathsf{R}$-matrix $\mathsf{S}$ of a square in the checkerboard chiral
Potts model,\cite{BPA} see Figs.~\ref{fig4} and \ref{fig7}, is the
intertwiner of two cyclic representations. Their procedure is sketched
in Fig.~\ref{fig8}.

They start from the six-vertex model $\mathsf{R}$-matrix at an $N$th
root of 1, with $N$ odd. The corresponding Yang--Baxter equation is
well-known with spin $\half$ highest-weight representations on all
legs. They then look for an $\mathsf{R}$-matrix $\mathsf{L}$
intertwining a highest-weight and a cyclic representation. This is
solved from a new Yang--Baxter equation that is quadratic in
$\mathsf{L}$.\footnote{Korepanov had earlier obtained the first
part of this construction, but his work has not yet been
published.\cite{Korv}}\ \ \ Next, the Yang--Baxter equation with one
highest-weight and two cyclic rapidity lines is a linear equation for
the intertwiner
$\mathsf{S}$ of two cyclic representations. This intertwiner
$\mathsf{S}$ satisfies a Yang--Baxter equation with only cyclic
representations on the legs. The original result\cite{BPA} for
$\mathsf{S}$ is obtained in a suitable gauge with a proper choice of
parameters.

The above illustrates the group theoretical significance of the chiral
Potts model. It gives a standard example of the affine quantum group
U$_q\widehat{\rm sl(2)}$ at root of unity in the minimal cyclic
representation.\cite{dCK} Intertwiners of certain bigger irreducible
cyclic representations are given by chiral Potts model partition
functions, as was found by Tarasov.\cite{Tara}

Smaller representations ({\it i.e.}\ semicyclic\cite{GRS}$^-$\cite{Schn}
and highest-weight) follow by the reduction of products of cyclic
representations in the construction of Baxter {\it et al}.\cite{BBP}
Elsewhere we shall present more details on how this relates several
root-of-unity models. Relating $p'$ in Fig.~\ref{fig5} by a special
automorphism with $p$, the product of transfer matrices splits
as\cite{BBP}
\begin{equation}
{\mathsf{T}}_q\,\hat{\mathsf{T}}_r=
B^{(j)}_{pp'q}\,{\mathsf{X}}^{-k}\,\tau^{(j)}(t_q)+
B^{(N-j)}_{p'pq}{\mathsf{X}}^l\,\tau^{(N-j)}(t_r).
\label{TT-split}
\end{equation}
Here the transfer matrix $\tau^{(j)}$ is made up of
$\mathsf{L}$-operators intertwining a cyclic and a spin
$s=\textstyle{j-1\over2}$ representation, the $B$'s are scalars,
$t_q\equiv x_qy_q$, and powers of the spin shift operator $\mathsf{X}$
come in depending on the automorphism.

Doing the same process $r=q'$ in the other direction, one obtains
several fusion relations for the $\tau^{(j)}$ transfer
matrices\cite{BBP}
\begin{eqnarray}
\tau^{(j)}(t_q)\,\tau^{(2)}(\omega^{j-1}t_q)&=&
z(\omega^{j-1}t_q)\,\tau^{(j-1)}(t_q)+\tau^{(j+1)}(t_q),\nonumber\\
\tau^{(j)}(\omega t_q)\,\tau^{(2)}(t_q)&=&
z(\omega t_q)\,\tau^{(j-1)}(\omega^2t_q)+\tau^{(j+1)}(t_q),
\label{tau-split}
\end{eqnarray}
where $z(t)$ is a known scalar function.

\section{Selected Exact Results}
\noindent
There are several exact results for the free energy of the integrable
chiral Potts model and we quote here a recent result of Baxter in the
scaling regime,\cite{B-scal}
\begin{eqnarray}
f-f_{\rm FZ}&\equiv&f-f_{\rm c}=
-{(N-1)k^2\over2N\pi}(u_q-u_p)\cos(u_p+u_q)\nonumber\\
&+&{k^2\over4\pi^2}\sin(u_q-u_p)
\sum_{j=1}^{<N/2}{\tan(\pi j/N)\over j}
B\bigg(1+{j\over N},{1\over2}\bigg)^2
\bigg({k\over2}\bigg)^{4j/N}\nonumber\\
&&+{\rm O}(k^4\log k),
\quad\hbox{if}\quad k^2\sim T_{\rm c}-T\to0,\quad\alpha=1-{2\over N}.
\label{scaling}
\end{eqnarray}
For the order parameters we have a general conjecture\cite{AMPT} in
the ordered state,
\begin{equation}
\langle\sigma_0^n\rangle = {(1-{k'}^2)}^{\beta_n},\quad
\beta_n={n(N-n)\over 2N^2},\quad
(1\le n\le N-1,\quad \sigma_0^N=1),
\label{eq:orderpar}
\end{equation}
which still remains to be proved.

Using Baxter's results\cite{Bax} we have a very
explicit formula\cite{AP-O} for the interfacial tensions,
\begin{equation}
{\epsilon_r\over k_{\rm B}T}=
{8\over\pi}\,\int_0^{\eta}dy\,{\sin(\pi r/N)\over
1+2y\cos(\pi r/ N)+y^2}\,
\hbox{artanh}\sqrt{\eta^N-y^N\over 1-\eta^N y^N},
\label{intf-T}
\end{equation}
in the fully symmetric, $W\equiv\wb$, integrable chiral Potts model,
{\it i.e.}\ diagonal chiral fields. Here
$\eta\equiv[(1-k')/(1+k')]^{1/N}$ is a temperature-like variable and
$r$ is the difference of the spin variables across the interface,
($r=1,\cdots,N-1$).

In the low-temperature region we have $k'\to0$ and
we can expand (\ref{intf-T}) as\cite{AP-O}
\begin{equation}
{\epsilon_r\over k_{\rm B}T}=
-{2r\over N}\,\log{k'\over 2}-
\log\biggl({1\over\cos^2 r\bar\lambda}
\prod_{j=1}^{[r/2]}\,{\cos^4 (r-2j+1)\bar\lambda
\over\cos^4 (r-2j)\bar\lambda}\biggr)
+{\rm O}(k'),
\label{intf-low}
\end{equation}
where $\bar\lambda\equiv \pi/(2N)$ and the constant term comes
from a dilogarithm integral. For the critical region,
$\eta\approx (k/2)^{2/N}\sim (T_{\rm c}-T)^{1/N}\to 0$,
(\ref{intf-T}) reduces to\cite{AP-O}
\begin{eqnarray}
{\epsilon_r\over k_{\rm B}T} & = &
{8\sin({\pi r/N})B({1/N},1/2)\over\pi(N+2)}\,\eta^{1+N/2} \nonumber \\
 & - &{8\sin({2\pi r/N})B({2/N},1/2)\over\pi(N+4)}\,\eta^{2+N/2}+
{\rm O}(\eta^{3+N/2}) \nonumber \\
 & \approx & \eta^{N\mu} D_r(\eta)
=\eta^{N\mu} D_r(\Delta/\eta^{N\phi}),
\label{intf-c}
\end{eqnarray}
where we have assumed the existence of scaling function $D_r$ depending
on $\eta$ and $\Delta\sim (T_{\rm c}-T)^{1/2}$, the chiral field
strength on the integrable line. From this, we have the critical
exponents
\begin{equation}
\mu={1\over 2}+{1\over N}=\nu, \quad \phi={1\over 2}-{1\over N}.
\label{intf-ex}
\end{equation}

We note that the above dilogarithm identity involves the dilogarithms
at $2N$-th roots of unity and has been discovered numerically first.
Apart from direct proofs for small $N=2,3,4$, only an indirect
proof$\,$\cite{AP-O} by the Bethe Ansatz exists.

\section{Basic Hypergeometric Series at Root of Unity}
\noindent
The basic hypergeometric hypergeometric series is
defined\cite{Bailey,Slater} as
\begin{equation}
\hypg{\alpha_1,\cdots,\alpha_{p+1}}
{\phantom{\alpha_1}\beta_1,\cdots,\beta_p}{z}=
\sum_{l=0}^{\infty}{{(\alpha_1;q)_l\cdots(\alpha_{p+1};q)_l}
\over{(\beta_1;q)_l\cdots(\beta_{p};q)_l(q;q)_l}}\,z^{l},
\label{eq:hypg}\end{equation}
where
\begin{equation}
(x;q)_l\equiv\cases{1,\quad l=0,\cr
(1-x)(1-xq)\cdots(1-xq^{l-1}),\quad l>0,\cr
1/[(1-xq^{-1})(1-xq^{-2})\cdots(1-xq^{l})],\quad l<0.}
\label{eq:qp}\end{equation}
Setting $\alpha_{p+1}=q^{1-N}$ and $q\to\omega\equiv e^{2\pi i/N}$,
we get
\begin{equation}
\hypg{\omega,\alpha_1,\cdots,\alpha_{p}}
{\phantom{\omega,}\beta_1,\cdots,\beta_p}{z}=
\sum_{l=0}^{N-1}
{{(\alpha_1;\omega)_l\cdots(\alpha_{p};\omega)_l}\over
{(\beta_1;\omega)_l\cdots(\beta_{p};\omega)_l}}\,z^{l}.
\label{eq:hypc}
\end{equation}

We note
\begin{equation}
(x;\omega)_{l+N}=(1-x^N)(x;\omega)_{l},\qquad\mbox{and}\quad
(\omega;\omega)_l=0,\quad l\ge N.
\label{eq:wp}
\end{equation}
So if we also require
\begin{equation}
z^N=\prod_{j=1}^p\gamma_j^N,\qquad
{\gamma_j}^N={1-\beta_j^N\over1-\alpha_j^N},
\label{eq:defg}\end{equation}
we obtain a ``cyclic basic hypergeometric function" with summand
periodic mod $N$.

For us the Saalsch\"utz case, defined by
\begin{equation}
z=q={\beta_1\cdots\beta_p\over\alpha_1\cdots\alpha_{p+1}}
\qquad\mbox{or}\quad
\omega^2\alpha_1\alpha_2\cdots\alpha_{p}=\beta_1\beta_2\cdots\beta_p,
\quad z=\omega,
\label{eq:saalc}
\end{equation}
is important, but details on this will be given elsewhere.

The theory of cyclic hypergeometric series is intimately related with
the theory of the integrable chiral Potts model and several identities
appear hidden in the literature. We note that our notations here, which
match up nicely with the classical definitions of basic hypergeometric
functions,\cite{Bailey,Slater} differ from those of Bazhanov\break {\it
et al}.,\cite{BB}$^-$\cite{SMS} who are using an upside-down version of
the $q$-Pochhammer symbol $(x;q)_l$. Hence, our definition of the
cyclic hypergeometric series differs from the one of Sergeev {\it et
al}.,\cite{SMS} who also use homogeneous rather than more compact affine
variables. So comparing our results with theirs is a little
cumbersome.

\pagebreak
\subsection{Integrable chiral Potts model weights}
\noindent
The weights of the integrable chiral Potts model can be written in
product form
\begin{equation}
{W(n)\over W(0)}=\gamma^{n}\,{(\alpha;\omega)_n\over(\beta;\omega)_n},
\qquad\gamma^N={1-\beta^N\over1-\alpha^N}.
\label{eq:w}
\end{equation}
This is periodic with period $N$ as follows from (\ref{eq:wp}).

The dual weights are given by Fourier transform,\cite{WuWang} {\it i.e.}
\begin{equation}
W^{({\rm f})}(k)=\sum_{n=0}^{N-1}\omega^{nk}\,W(n)=
\hypp{1}{\alpha}{\beta}{\gamma\,\omega^k}\,W(0).
\label{eq:wf}
\end{equation}
Using the recursion formula
\begin{equation}
W(n)\,(1-\beta\,\omega^{n-1})=W(n-1)\,\gamma\,(1-\alpha\,\omega^{n-1})
\label{eq:ww}\end{equation}
and its Fourier transform, we find
\begin{equation}
{W^{({\rm f})}(k)\over W^{({\rm f})}(0)}=
{\displaystyle{\hypp{1}{\alpha}{\beta}{\gamma\,\omega^k}}\over
\displaystyle{\hypp{1}{\alpha}{\beta}{\gamma}}}
=\left({\omega\over\beta}\right)^k
{(\gamma;\omega)_k\over(\omega\alpha\gamma/\beta;\omega)_k}.
\label{eq:rec1}\end{equation}
This relation is equivalent to the one originally found in
1987,\cite{BPA,AP-tani} It shows that dual weights also satisfy
(\ref{eq:w}).

With one more Fourier transform we get
\begin{equation}
{N\over\displaystyle{\hypp{1}{\alpha}{\beta}{\gamma}}}=
{\hypp{1}{\phantom{q/b}\gamma}{\omega\alpha\gamma/\beta}
{{\omega\over\beta}}}
={\hypp{1}{\beta/\alpha\gamma}{\omega/\gamma}{\alpha}}.
\label{eq:wff}
\end{equation}
Also, we can show that
\begin{equation}
{\displaystyle{\hypp{1}{\alpha\,\omega^m}{\beta\,\omega^n}
{\gamma\,\omega^k}}\over
\displaystyle{\hypp{1}{\alpha}{\beta}{\gamma}}}=
{(\omega/\beta)^{k}(\beta;\omega)_n(\gamma;\omega)_k
(\omega\alpha/\beta;\omega)_{m-n}
\over{(\gamma\,\omega^k)}^{n}(\alpha;\omega)_m
(\omega\alpha\gamma/\beta;\omega)_{m-n+k}}.
\label{eq:rec4}\end{equation}
This equation has been proved by Kashaev {\it et al}.\cite{KMS} in
other notation and is valid for all values of the arguments,
provided condition (\ref{eq:w}) holds.

\subsection{Baxter's summation formula}
\noindent
{}From Baxter's work\cite{B-399} we can infer the identity
\begin{equation}\begin{array}{ll}
\displaystyle{\hypp{1}{\alpha}{\beta}{\gamma}}&=
\displaystyle{\Phi_0\,\sqrt{N}\,(\omega/\beta)^{{1\over2}(N-1)}}\cr
&\displaystyle{\times\,
\prod_{j=1}^{N-1}\left[{(1-\omega^{j+1}\alpha/\beta)(1-\omega^j\gamma)
\over(1-\omega^j\alpha)(1-\omega^{j+1}/\beta)
(1-\omega^{j+1}\alpha\gamma/\beta)}\right]^{j/N}\!\!,}\cr
\end{array}\label{eq:hyp1}
\end{equation}
valid up to an $N$-th root of unity, while
\begin{equation}
\Phi_0\equiv e^{i\pi(N-1)(N-2)/12N},\qquad
\gamma^N={1-\beta^N\over1-\alpha^N}.
\label{eq:phi0}
\end{equation}
Introducing a function
\begin{equation}
p(\alpha)=\prod_{j=1}^{N-1}(1-\omega^{j}\alpha)^{j/N},\qquad p(0)=1,
\label{eq:defp}\end{equation}
we can rewrite the identity as
\begin{equation}
\hypp{1}{\alpha}{\beta}{\gamma}=\omega^d N^{1\over2}\Phi_0
\left({\omega\over\beta}\right)^{{1\over2}(N-1)}
{p(\omega \alpha/\beta)p(\gamma)
\over p(\alpha)p(\omega/\beta)
p(\omega\alpha\gamma/\beta)},
\label{eq:hyp1d}
\end{equation}
with $\omega^d$ determined by the choice of the branches. The LHS of
(\ref{eq:hyp1d}) is single valued in $\alpha$, $\beta$, and $\gamma$,
whereas the RHS has branch cuts. It is possible to give a precise
prescription for $d$, assuming that $p(\alpha)$ has branch cuts along
$[\omega^j,\infty)$ for $j=1,\ldots,N-1$, following straight lines
through the origin; but these details will be presented elsewhere.

Using this relation (\ref{eq:hyp1d}) and classical identities we can get
many relations for ${}_3\Phi_2,\,{}_4\Phi_3,\,{}_5\Phi_4,\ldots$,
including the star--triangle equation and the tetrahedron equation of
the Bazhanov--Baxter model. Some of these results have been given very
recently in different notation.\cite{SMS}

\subsection{Outline of proof of Baxter's identity}
\noindent
As (\ref{eq:hyp1d}) is crucial in the theory, we briefly sketch a
proof here. Also we note that in this subsection, we choose the
normalization $W(0)=1$ for the weights given by (\ref{eq:w}), rather
than $\prod W(j)=1$.

We note that the $W^{({\rm f})}(k)$ are eigenvalues of the cyclic
$N\times N$ matrix
\begin{equation}
{\mathsf{M}}\equiv\left(\matrix{W(0)&W(1)&\ldots&W(N-1)\cr
W(-1)&W(0)&\ldots&W(N-2)\cr.
\vdots&\vdots&\ddots&\vdots\cr
W(1-N)&W(2-N)&\ldots&W(0)\cr}\right),
\label{eq:matr}
\end{equation}
and
\begin{equation}
\det{\mathsf{M}}=
\prod_{j=0}^{N-1}W^{(f)}(j)=[W(0)^{(f)}]^N\prod_{j=1}^{N-1}
\left[{W^{(f)}(j)\over W^{(f)}(0)}\right].
\label{eq:det}
\end{equation}
So if we can calculate $\det\mathsf{M}$ directly, we obtain a result
for $[W(0)^{(f)}]^N$, giving us the proof of the identity.

Using
\begin{equation}
W(n)=\gamma^n\,\prod_{j=0}^{n-1}{1-\alpha\omega^j\over1-\beta\omega^j},
\quad W(N-n)=W(-n)=
\gamma^{-n}\,\prod_{j=-n}^{-1}{1-\beta\omega^j\over1-\alpha\omega^j},
\label{eq:detsub}
\end{equation}
we can rewrite
\begin{equation}
\det{\mathsf{M}}=\prod_{l=0}^{N-1}\Biggl[\prod_{j=0}^{N-2-l}
(1-\omega^{j}\beta)^{-1}
\prod_{j=-l}^{-1}(1-\omega^{j}\alpha)^{-1}\Biggr]\,
\det{\mathsf{E}}^{(0)},
\label{eq:d}\end{equation}
where the elements of ${\mathsf{E}}^{(0)}$ are polynomials, so that
$\det{\mathsf{E}}^{(0)}$ is also a polynomial. We can define more
general matrix elements ${\mathsf{E}}^{(m)}$ for $m=0,\dots,N-1$ by
\begin{equation}
E^{(m)}_{k,l}=\prod_{j=-k+1}^{l-k-1}
(1-\omega^{j}\alpha)
\prod_{j=l+m-k}^{N-k-1}(1-\omega^{j}\beta),
\label{eq:ed}
\end{equation}
which satisfy the recursion relation
\begin{equation}
E^{(m)}_{k,l}-E^{(m)}_{k,l+1}=\omega^{l-k}\,(\alpha-\omega^m\beta)\,
E^{(m+1)}_{k,l}.
\label{eq:em}
\end{equation}
Subtracting the pairs of consecutive columns of $\det{\mathsf{E}}^{(0)}$
in (\ref{eq:d}), and using (\ref{eq:em}), we can pull out some of the
zeros leaving a determinant with $N-1$ columns from ${\mathsf{E}}^{(1)}$
and the last column from ${\mathsf{E}}^{(0)}$. Repeating this process,
we arrive at
\begin{equation}
\det{\mathsf{E}}^{(0)}=
\prod_{m=0}^{N-2}(\alpha-\omega^{m}\beta)^{N-1-m}
\cdot\det{\mathsf{F}}.
\label{eq:dete2}
\end{equation}
Here the matrix $\mathsf{F}$ is defined such that its $j$th column is
the $j$th column of matrix ${\mathsf{E}}^{(N-j)}$.

{}From a simple polynomial degree count we conclude that
$\det{\mathsf{F}}$ has to be a constant. Noting that
${\mathsf{E}}^{(0)}$ is triangular in the limit $\alpha\to1$, we find
\begin{eqnarray}
\det{\mathsf{F}}=\prod_{j=1}^{N-1}(1-\omega^{j})^j=
\Phi_0^N\,N^{{1\over2}N},
\label{eq:fv}
\end{eqnarray}
in which $\Phi_0$ is given by (\ref{eq:phi0}). Hence, we can complete
the proof of the identity.

\subsection{Further identities}
\noindent
Several other identities can be derived using the above identity and
the classical Jackson identity.\cite{Bailey,Slater} One thus generates
the fundamental identities for the weights of the Baxter--Bazhanov model
and the sl$(n)$ chiral Potts model. More precisely, the Boltzmann
weights of a cube in the Baxter--Bazhanov model are proportional to
${}_3\Phi_2$'s, so all identities for the weights of this model are
also identities for cyclic ${}_3\Phi_2$'s.

Without any restriction on the parameters, we can derive
\begin{equation}
\hypp{2}{\alpha_1,\alpha_2}{\beta_1,\beta_2}{z}=N^{-1}
\sum_{k=0}^{N-1}\hypp{1}{\alpha_1}{\beta_1}{\omega^{-k}\gamma_1}
\hypp{1}{\alpha_2}{\beta_2}{{\omega^k \gamma_2}},
\label{eq:hyp2}\end{equation}
where $z=\gamma_1\gamma_2$ and $\gamma_i$ is defined in (\ref{eq:defg}).
We have only used the convolution theorem so far. Next, we can use
(\ref{eq:rec1}) so that we can perform the sum in (\ref{eq:hyp2}). This
way we obtain the transformation formula
\begin{equation}
\hypp{2}{\alpha_1,\alpha_2}{\beta_1,\beta_2}{z}=A\,\,\,
\hypp{2}{{z/\gamma_1},
{\hphantom{\omega}{\beta_1/\alpha_1\gamma_1}\hphantom{\omega}}}
{{\omega/\gamma_1},{\omega\alpha_2z/\beta_2\gamma_1}}
{{\omega\alpha_1\over\beta_2}},
\label{eq:h2t1}
\end{equation}
where the constant $A$ can be written in several different forms,
either with ${}_2\Phi_1$'s or with $p(x)$'s using (\ref{eq:hyp1d}).
{}From this identity one can generate the symmetry relations of
the cube in the Baxter--Bazhanov model under the 48 elements of the
symmetry group of the cube, see also the recent work of Sergeev
{\it et al}.\cite{SMS}

In addition, one can work out relations for ${}_4\Phi_3$ and higher.
One of these many\break relations, {\it i.e.}\ a Saalsch\"utzian
${}_4\Phi_3$ identity, is the star--triangle equation of the\break
integrable chiral Potts model,\cite{BPA} or its Fourier
transform\cite{AMPTY}
\begin{equation}
V_{prq}(a,b;n)\,\wb_{qr}^{({\rm f})}(n)=
R_{pqr}\,V_{pqr}(a,b;n)\,W_{qr}(a-b).
\label{eq:STEV}
\end{equation}

More detail will be presented elsewhere. We also have to refer the
reader to the recent work of Stroganov's group\cite{SMS} which uses
fairly different notations in their appendix. Their higher identities
also involve Saalsch\"utzian cyclic hypergeometric functions, albeit
that that is hard to recognize.

As a conclusion, we may safely state that the existence of all these
cyclic\break hypergeometric identities is the mathematical reason
behind the integrable chiral Potts family of models.

\nonumsection{Acknowledgements}
\noindent
We thank Professor F.\ Y.\ Wu for his hospitality and Professors
R.\ J.\ Baxter and Yu.\ G.\ Stroganov for much advice.
This work has been supported in part by NSF
Grants No.\ PHY-93-07816 and PHY-95-07769.

\nonumsection{References}


\noindent
\begin{thebibliography}{99}

\bibitem{WuWang}F.\ Y.\ Wu and Y.\ K.\ Wang,
\Journal{\myem J.\ Math.\ Phys.}{17}{439}{1976}.

\bibitem{Os}S.\ Ostlund,
\Journal{\myem Phys.\ Rev.}{B24}{398}{1981}.

\bibitem{Hu}D.\ A.\ Huse,
\Journal{\myem Phys.\ Rev.}{B24}{5180}{1981}.

\bibitem{HSF}D.\ A.\ Huse, A.\ M.\ Szpilka, and M.\ E.\ Fisher,
\Journal{\myem Physica}{A121}{363}{1983}.

\bibitem{dN}M.\ P.\ M.\ den Nijs,
in {\myem Phase Transitions and Critical Phenomena},
eds.\ C.\ Domb and\break J.\ L.\ Lebowitz
(Academic, London, 1988), Vol.\ 12, p.\ 219.

\bibitem{Selke}W.\ Selke,
\Journal{\myem Phys.\ Repts.}{170}{213}{1988}.

\bibitem{AP-O}H.\ Au-Yang and J.\ H.\ H.\ Perk,
\Journal{\myem J.\ Stat.\ Phys.}{78}{17}{1995}.

\bibitem{Barber}M.\ N.\ Barber,
\Journal{\myem J.\ Phys.}{A15}{915}{1982}.

\bibitem{HF}D.\ A.\ Huse and M.\ E.\ Fisher,
\Journal{\myem J.\ Phys.}{C15}{L585}{1982}.

\bibitem{AP-SP}H.\ Au-Yang and J.\ H.\ H.\ Perk,
in {\myem STATPHYS 19}, ed.\ Hao Bailin
(World Scientific, Singapore, 1996), p.\ 421.

\bibitem{AP-F}H.\ Au-Yang and J.\ H.\ H.\ Perk,
\Journal{\myem Physica}{A177}{139}{1991}.

\bibitem{AMPTY}H.\ Au-Yang, B.\ M.\ McCoy, J.\ H.\ H.\ Perk,
S.\ Tang, and M.\ L.\ Yan,
\Journal{\myem Phys.\ Lett.}{A123}{219}{1987}.

\bibitem{AMPT-S}H.\ Au-Yang, B.\ M.\ McCoy, J.\ H.\ H.\ Perk,
and S.\ Tang,
in {\myem Algebraic Analysis}, Vol.\ 1,
eds.\ M. Kashiwara and T. Kawai
(Academic, San Diego, 1988), p.\ 29.

\bibitem{BPA}R.\ J.\ Baxter, J.\ H.\ H.\ Perk, and H.\ Au-Yang,
\Journal{\myem Phys.\ Lett.}{A128}{138}{1988}.

\bibitem{AP-tani}H.\ Au-Yang and J.\ H.\ H.\ Perk,
\Journal{\myem Adv.\ Stud.\ Pure Math.}{19}{57}{1989}.

\bibitem{AMP}G.\ Albertini, B.\ M.\ McCoy, and J.\ H.\ H.\ Perk,
\Journal{\myem Adv.\ Stud.\ Pure Math.}{19}{1}{1989};
\Journal{\myem Phys.\ Lett.}{A135}{159}{1989},
\Journal{$\!\!$}{A139}{204}{1989}.

\bibitem{FZ}V.\ A.\ Fateev and A.\ B.\ Zamolodchikov,
\Journal{\myem Phys.\ Lett.}{A92}{37}{1982}.

\bibitem{AMPT}G.\ Albertini, B.\ M.\ McCoy, J.\ H.\ H.\ Perk,
and S.\ Tang,
\Journal{\myem Nucl.\ Phys.}{B314}{741}{1989}.

\bibitem{MR}B.\ M.\ McCoy and S.-S.\ Roan,
\Journal{\myem Phys.\ Lett.}{A150}{347}{1990}.

\bibitem{B-399}R.\ J.\ Baxter,
\Journal{\myem J.\ Stat.\ Phys.}{52}{639}{1988}.

\bibitem{B-sup}R.\ J.\ Baxter,
\Journal{\myem J.\ Stat.\ Phys.}{57}{1}{1989}.

\bibitem{B-eiva}R.\ J.\ Baxter,
in {\myem Calculation of the Eigenvalues of the Transfer Matrix of the
Chiral Potts Model,}
Proc.\ of Fourth Asia Pacific Physics Conference
(World Scientific, Singapore, 1991), p.\ 42.

\bibitem{B-scal}R.\ J.\ Baxter,
\Journal{\myem J.\ Stat.\ Phys.}{82}{1219}{1996}.

\bibitem{BBP}R.\ J.\ Baxter, V.\ V.\ Bazhanov, and J.\ H.\ H.\ Perk,
\Journal{\myem Intern.\ J.\ Mod.\ Phys.}{B4}{803}{1990}.

\bibitem{Bax}R.\ J.\ Baxter,
\Journal{\myem J.\ Stat.\ Phys.}{73}{461}{1993},
\Journal{\myem J.\ Phys.}{A27}{1837}{1994}.

\bibitem{oRB}R.\ J.\ Baxter,
M.\ J.\ O'Rourke and R.\ J.\ Baxter,
\Journal{\myem J.\ Stat.\ Phys.}{82}{1}{1996}.

\bibitem{AP-multi}H.\ Au-Yang and J.\ H.\ H.\ Perk,
\Journal{\myem Intern.\ J.\ Mod.\ Phys.}{A7}{Suppl.\ 1B,
1007, 1025}{1992}.

\bibitem{MPS}W.\ S.\ McCullough, J.\ H.\ H.\ Perk, and H.\ L.\ Scott,
\Journal{\myem J.\ Chem.\ Phys.}{93}{6070}{1990}.

\bibitem{AP-inf}H.\ Au-Yang and J.\ H.\ H.\ Perk,
\Journal{\myem Intern.\ J.\ Mod.\ Phys.\
(Proc.\ Suppl.)}{3A}{430}{1993}.

\bibitem{dCK}C.\ De Concini and V.\ G.\ Kac,
in {\myem Operator Algebras, Unitary Representations,\break Enveloping
Algebras, and Invariant Theory},
eds.\  A.\ Connes, M.\ Duflo, A.\ Joseph, and R.\ Rentschler
(Birkh\"auser, Boston, 1990), p.\ 471.

\bibitem{HKdN}S.\ Howes, L.\ P.\ Kadanoff, and M.\ den Nijs,
\Journal{\myem Nucl.\ Phys.}{B215[FS7]}{169}{1983}.

\bibitem{vGR}G.\ von Gehlen and V.\ Rittenberg,
\Journal{\myem Nucl.\ Phys.}{B257[FS14]}{351}{1985}.

\bibitem{PD}J.\ H.\ H.\ Perk,
\Journal{\myem Proc.\ Symp.\ Pure Math.}{49}{part 1, 341}{1989}.

\bibitem{Dav}B.\ Davies,
\Journal{\myem J.\ Math.\ Phys.}{32}{2945}{1991}.

\bibitem{BB}V.\ V.\ Bazhanov and R.\ J.\ Baxter,
\Journal{\myem J.\ Stat.\ Phys.}{69}{453}{1992},
\Journal{$\!\!$}{71}{839}{1993}.

\bibitem{KMS}R.\ M.\ Kashaev, V.\ V.\ Mangazeev,
and Yu.\ G.\ Stroganov,
\Journal{\myem Intern.\ J.\ Mod.\ Phys.}{A8}{587,1399}{1993}.

\bibitem{MSS}V.\ V.\ Mangazeev, S.\ M.\ Sergeev, and Yu.\ G.\ Stroganov,
\Journal{\myem Intern.\ J.\ Mod.\ Phys.}{A9}{5517}{1994};
\Journal{\myem Mod.\ Phys.\ Lett.}{A10}{279}{1995}.

\bibitem{BMS}G.\ E.\ Boos, V.\ V.\ Mangazeev, and S.\ M.\ Sergeev,
\Journal{\myem Intern.\ J.\ Mod.\ Phys.}{A10}{4041}{1996}.

\bibitem{SBMS}S.\ M.\ Sergeev, G.\ E.\ Boos, V.\ V.\ Mangazeev,
and Yu.\ G.\ Stroganov,
\Journal{\myem Mod.\ Phys.\ Lett.}{A11}{491}{1996}.

\bibitem{Bo}G.\ E.\ Boos,
\Journal{\myem Intern.\ J.\ Mod.\ Phys.}{A11}{313}{1996}.

\bibitem{Huzn}Z.\ N.\ Hu,
\Journal{\myem Intern.\ J.\ Mod.\ Phys.}{A9}{5201}{1994};
\Journal{\myem Phys.\ Lett.}{A197}{387}{1995}.

\bibitem{HuH}Z.\ N.\ Hu and B.\ Y.\ Hou,
\Journal{\myem J.\ Stat.\ Phys.}{79}{759}{1995},
\Journal{$\!\!$}{82}{633}{1996}.

\bibitem{SMS}S.\ M.\ Sergeev, V.\ V.\ Mangazeev, and Yu.\ G.\ Stroganov,
\Journal{\myem J.\ Stat.\ Phys.}{82}{31}{1996}.

\bibitem{DJMM}E.\ Date, M.\ Jimbo, K.\ Miki and T.\ Miwa,
\Journal{\myem Comm.\ Math.\ Phys.}{137}{133}{1991};
\Journal{\myem Publ.\ RIMS, Kyoto Univ.}{27}{347, 639}{1991}.

\bibitem{BKMS}V.\ V.\ Bazhanov, R.\ M.\ Kashaev, V.\ V.\ Mangazeev,
and Yu.\ G.\ Stroganov,
\Journal{\myem Comm.\ Math.\ Phys.}{138}{393}{1991}.

\bibitem{KMN}R.\ M.\ Kashaev, V.\ V.\ Mangazeev and T.\ Nakanishi,
\Journal{\myem Nucl.\ Phys.}{B362}{563}{1991}.

\bibitem{KM}R.\ M.\ Kashaev and V.\ V.\ Mangazeev,
\Journal{\myem Mod.\ Phys.\ Lett.}{A7}{2827}{1992}.

\bibitem{AP-C}H.\ Au-Yang and J.\ H.\ H.\ Perk,
\Journal{\myem Physica}{A228}{78}{1996}. 

\bibitem{BS}V.\ V.\ Bazhanov and Yu.\ G.\ Stroganov,
\Journal{\myem J.\ Stat.\ Phys.}{59}{799}{1990}.

\bibitem{Korv}I.\ G.\ Korepanov,
{\myem unpublished reports} (in Russian, for more details see also
preprint hep-th/941006).

\bibitem{Tara}V.\ O.\ Tarasov,
\Journal{\myem Int.\ J.\ Mod.\ Phys.}{A7}{Suppl.\ 1B, 963}{1992};
\Journal{\myem Comm.\ Math.\ Phys.}{158}{459}{1993}.

\bibitem{GRS}C.\ G\'omez, M.\ Ruiz-Altaba and G.\ Sierra,
\Journal{\myem Phys.\ Lett.}{B265}{95}{1991}.

\bibitem{IU}I.\ T.\ Ivanov and D.\ B.\ Uglov,
\Journal{\myem Phys.\ Lett.}{A167}{459}{1992}.

\bibitem{Schn}A.\ Schnizer,
\Journal{\myem J.\ Math.\ Phys.}{34}{4340}{1993};
\Journal{\myem Comm.\ Math.\ Phys.}{163}{293}{1994}.

\bibitem{Bailey}W.\ N.\ Bailey,
{\myem Generalized Hypergeometric Series}
(Cambridge Univ.\ Press, Cambridge, 1935), pp.\ 65--70.

\bibitem{Slater}L.\ J.\ Slater,
{\myem Generalized Hypergeometric Functions}
(Cambridge Univ.\ Press, Cambridge, 1966), Ch.\ 3, App.\ II, IV.

\end{thebibliography}
\end{document}